\documentclass[oneside,reqno,12pt]{amsart}
\usepackage[greek,english]{babel}
\usepackage{amsthm}
\usepackage{amsbsy}
\usepackage{amsfonts}
\usepackage{graphicx}
\usepackage{hyperref}
\hypersetup{colorlinks=true, citecolor=blue, linkcolor=red}

\newtheorem{thm}{Theorem}

\newtheorem{rem}{Remark}

\numberwithin{equation}{section} \numberwithin{lem}{section}
\numberwithin{thm}{section} \numberwithin{cor}{section}
\numberwithin{pro}{section} \numberwithin{rem}{section}

\begin{document}
\title[A one-dimensional symmetry  result for the Fisher-KPP equation]{A one-dimensional symmetry  result for entire solutions to the Fisher-KPP equation}
\author{Christos Sourdis}
\address{National and Kapodistrian University of Athens, Department of Mathematics, Athens, Greece.
}
\email{sourdis@uoc.gr}

\date{\today}
\begin{abstract}
 We consider the  Fisher-KPP reaction-diffusion equation in the whole space.
 We prove that if a solution has, to main order and for all times (positive and negative), the same exponential decay as a planar traveling wave with speed larger than the minimal one  at its leading edge, then   it has to coincide with the aforementioned traveling wave.
 \end{abstract}
 \maketitle

\section{Introduction}
The Fisher-KPP equation
\begin{equation}\label{eqEq}
u_{t}=\Delta u +f(u),\ x\in \mathbb{R}^N,\ t\in \mathbb{R}, \ N\geq 1,
\end{equation}
with $f(u)=u(1-u)$
appears in the context of population dynamics to describe the spatial spread of an advantageous allele (see \cite{F,KPP}).
For every wave speed $c\geq 2$ and  $\eta\in \mathbb{S}^{N-1}$  it admits a unique, up to translations, planar \emph{traveling wave solution} of the form
\[ u(x,t)=V(\xi)\ \textrm{with}\ \xi=\eta \cdot x-ct\]
such that
\[V_\xi<0,\ \lim _{\xi \to -\infty }V(\xi)=1,\ \lim _{\xi \to +\infty }V(\xi)=0.\]
We point out that $V$ satisfies
\[
  v_{\xi\xi}+cv_{\xi}+f(v)=0,\ \xi\in \mathbb{R}.
\]
We note that such fronts move in the direction of $\eta$ with constant speed $c$, and their tail as $\xi\to +\infty$ is frequently referred to as the \emph{leading edge} of the wave.
More generally, there exists a $c_*>0$ such that, for each $c\geq c_*$, a completely analogous existence-uniqueness result holds  for $f\in C^1(\mathbb{R})$ satisfying
\begin{equation}\label{eqf}
  f(0)=f(1)=0,\ f>0\ \textrm{in}\ (0,1).
\end{equation}

In fact, if one further assumes that \begin{equation}\label{eqKPP}
                                       f(s)<f'(0)s,\ s\in (0,1),
                                     \end{equation} then $c_*=2\sqrt{f'(0)}$ holds and, possibly after a translation, $V$ satisfies the following asymptotic behaviour at its leading edge:
\begin{equation}\label{eqnara}
    V(\xi)=\left\{\begin{array}{ll}
                  \xi e^{r\xi}\left(1+o(1)\right) & \textrm{if}\ c=2\sqrt{f'(0)}, \\
                   &  \\
                  e^{r\xi}\left(1+o(1)\right) & \textrm{if}\ c>2\sqrt{f'(0)},
                \end{array}\right.
  \ \textrm{as}\ \xi\to +\infty,
\end{equation}
 where \begin{equation}\label{eqr} r=\frac{-c+\sqrt{c^2-4f'(0)}}{2}<0.\end{equation}
For the above properties, we refer to the introduction of \cite{nadiras} and the many references therein.

The purpose of this note  is to prove the following one-dimensional symmetry result.
\begin{thm}\label{thm1}
  We assume that  $f\in C^2(\mathbb{R})$ satisfies (\ref{eqf}) and  \begin{equation}\label{eqarxidia}
                                                                      f'(s)\leq f'(0),\ s\in (0,1).
                                                                    \end{equation} Let  $u$ be a solution to (\ref{eqEq}) with values in $(0,1)$ such that
for some  $\eta\in \mathbb{S}^{N-1}$ and $c>2\sqrt{f'(0)}$ it satisfies \begin{equation}\label{eqv}
  u(x,t)=
                  e^{r(\eta \cdot x-ct)}\left(1+o(1)\right)
    \ \textrm{as}\ \ \eta \cdot x-ct\to +\infty,
\end{equation}
 uniformly in $t\in \mathbb{R}$ and in the subspace of $\mathbb{R}^N$ that is orthogonal to $\eta$. Then,  we have
 \[
 u(x,t)\equiv V(\eta \cdot x-ct).
 \]
% with $y$ being the projection of $x$ to the orthogonal subspace to $\eta$, for some solution $U$ of
 %\[
 %\Delta u+c \eta \cdot \nabla u+f(u)=0\ \textrm{in}\ \mathbb{R}^N.
 %\]
 %In particular, if $N=1$ this implies that $u$ coincides with the planar traveling wave $V(\eta \cdot x-ct)$.
\end{thm}

\begin{rem}
  The assumptions on $f$ in Theorem \ref{thm1} imply that (\ref{eqKPP}) holds, and thus $c_*=2\sqrt{f'(0)}$.
\end{rem}
%\begin{rem}
%  As will be apparent from the proof, our arguments go through under the weaker assumption that $f\in C^2(\mathbb{R})$ satisfies (\ref{eqKPP}) amd $f'(0)\geq f'(s)$, $s\in (0,1)$ (see in particular (\ref{eqsuperaras}) below).
%\end{rem}

If $f$ is of class $C^1$, satisfies (\ref{eqf}) and $f'(0)>0$, $f'(1)<0$, it was shown in \cite[Thm. 3.5]{berestyckia} that if a solution of (\ref{eqEq}) satisfies
\begin{equation}\label{eqass}
  V(\eta\cdot x-ct)\leq u \leq V(\eta\cdot x-ct-a), \ x\in \mathbb{R}^N, \ t\in \mathbb{R},
\end{equation}for some $\eta \in \mathbb{S}^{N-1}$, $c\geq c_*$ and $a>0$,
  then \[u\equiv V(\eta\cdot x-ct-b)\ \textrm{for some}\ b\in [0,a].\]
If one further assumes that $f$ is concave in $[0, 1]$ and
$c > 2\sqrt{f'(0)}$ (recall that this is $c_*$ in the case of (\ref{eqKPP})), then the assertion (\ref{eqass}) follows from the weaker assumptions that
\[
 0 < u < 1\ \textrm{and}\ u(x,t) \to 1\ (\textrm{resp.}\ 0)\ \textrm{unif. as}\ \eta \cdot x  - ct \to -\infty \ (\textrm{resp.}\  +\infty),
\]
see \cite{nadiras} and \cite[Rem. 3.6]{berestyckia}. We refer to the former reference for the existence of an infinite-dimensional manifold of solutions to (\ref{eqEq}) which are not traveling waves (still for concave $f$).

Our method of proof is similar in spirit to the aforementioned references, in the sense that it relies on a sweeping argument and the strong maximum principle. However, there are substantial differences in the implementation of this general method.  Most notably, we apply our sweeping argument to the linearized equation of (\ref{eqEq})
(after switching to traveling wave coordinates). A main observation is that, thanks to (\ref{eqarxidia}), $e^{\eta \cdot x-ct}$ is a positive supersolution of the aforementioned equation if $c>2\sqrt{f'(0)}$.
In fact, we will sweep with this function.
If $c=2\sqrt{f'(0)}$ then $(\eta \cdot x-ct)e^{\eta \cdot x-ct}$ (recall (\ref{eqnara}))  is still a supersolution but it is sign changing, which is the reason why our proof breaks down in that case.
Related ideas in a different context can be found in our recent paper \cite{sourdaras}.

The rest of the paper is devoted to the proof of Theorem \ref{thm1}.
\section{Proof of Theorem \ref{thm1}}
\begin{proof}
Without loss of generality, we may assume that $\eta=(1,0,\cdots,0)$.

\underline{\textbf{Traveling wave coordinates.}}
It is natural to study $u$ in  traveling wave coordinates. For this purpose, with a slight abuse of notation, we write
\[
u(x,t)=u(x_1,\cdots,x_N,t)=u(\xi,y,t)\ \textrm{with} \ \xi=x_1-ct\ \textrm{and}\ y=(x_2,\cdots,x_N).
\]
In this frame of reference, $u$ solves
\begin{equation}\label{equPDE}
  u_t=u_{\xi\xi}+cu_{\xi}+\Delta_yu+f(u),\ \ \xi\in \mathbb{R}, \ y\in \mathbb{R}^{N-1},\ t\in \mathbb{R},
\end{equation}
and satisfies
\begin{equation}\label{equAS}
  u(\xi,y,t)= e^{r\xi}+w(\xi,y,t),
\end{equation}
where $w$ is such that
\begin{equation}\label{eqw}
w(\xi,y,t)=o( e^{r\xi})\ \textrm{as} \ \xi \to +\infty,
\end{equation}
uniformly in $y\in \mathbb{R}^{N-1}$ and $t\in \mathbb{R}$.
 In these coordinates, which we will use throughout the rest of the proof,  the assertion of the theorem reduces to \begin{equation}\label{eqdesire}
                                                                                                                      u_t\equiv 0\ \textrm{and} \ \nabla_yu\equiv 0.
                                                                                                                    \end{equation} We will only show the first identity of the above relation since the other one can be established in a completely analogous fashion.

\underline{\textbf{Gradient estimates.}}
Our next objective is to see what the asymptotic behaviour (\ref{eqv}) implies for $u_t$.
Since $ e^{r\xi}$ solves the linearized problem
\begin{equation}\label{eqSupreme}
v_{\xi\xi}+cv_{\xi}+f'(0)v=0,\ \xi\in \mathbb{R},
\end{equation}
we find that $w$ solves
\begin{equation}\label{eqwEq}
w_t-w_{\xi\xi}-cw_{\xi}-\Delta_yw=f( e^{r\xi}+w)-f'(0) e^{r\xi}.
\end{equation}
The righthand side of the above equation can be written as
\begin{equation}\label{eqmount}
f( e^{r\xi}+w)-f( e^{r\xi})+f( e^{r\xi})-f(0)-f'(0) e^{r\xi}\stackrel{(\ref{eqw})}{=}o( e^{r\xi})\ \textrm{as} \ \xi \to +\infty,
\end{equation}
uniformly in $y\in \mathbb{R}^{N-1}$ and $t\in \mathbb{R}$.
Then, by applying standard interior parabolic $W_p^{2,1}$ estimates (see for instance \cite[Thm. 7.22]{lieber}) in cylinders of the form $\mathcal{C}_{(\xi,y,t)}=\left\{|\Xi-\xi|<1,\ |Y-y|<1,\ |T-t|<1 \right\}$,
and  using  the
parabolic Sobolev embedding (see \cite[pgs. 80, 342]{lady}), we infer that
\begin{equation}\label{eqwGrad}
|\nabla_{\xi,y}w|=o( e^{r\xi})\ \textrm{as} \ \xi \to +\infty,
\end{equation}
uniformly in $y\in \mathbb{R}^{N-1}$ and $t\in \mathbb{R}$.

Let \begin{equation}\label{eqz}z=w_\xi.\end{equation}
Differentiation of (\ref{eqwEq}) with respect to $\xi$ yields
\[
z_t-z_{\xi\xi}-cz_{\xi}-\Delta_yz=f'( e^{r\xi}+w)(r e^{r\xi}+z)-f'(0)r e^{r\xi}\stackrel{(\ref{eqw}), (\ref{eqwGrad})}{=}o( e^{r\xi}),
\]
as $\xi \to +\infty$, uniformly in $y\in \mathbb{R}^{N-1}$ and $t\in \mathbb{R}$.
By the same procedure as before, we get
\begin{equation}\label{eqwGrad1}
|\nabla_{\xi,y}z|=o( e^{r\xi})\ \textrm{as} \ \xi \to +\infty,
\end{equation}
uniformly in $y\in \mathbb{R}^{N-1}$ and $t\in \mathbb{R}$.
For $i=2,\cdots,N$, let
\begin{equation}\label{eqomega}\omega=w_{x_i}.\end{equation}
Differentiation of (\ref{eqwEq}) now with respect to $x_i$ yields
\[
\omega_t-\omega_{\xi\xi}-c\omega_{\xi}-\Delta_y\omega=f'( e^{r\xi}+w)\omega\stackrel{(\ref{eqwGrad})}{=}o( e^{r\xi}),
\]
as $\xi \to +\infty$, uniformly in $y\in \mathbb{R}^{N-1}$ and $t\in \mathbb{R}$.
By working in the usual way, we obtain
\begin{equation}\label{eqwGrad2}
|\nabla_{\xi,y}\omega|=o( e^{r\xi})\ \textrm{as} \ \xi \to +\infty,
\end{equation}
uniformly in $y\in \mathbb{R}^{N-1}$ and $t\in \mathbb{R}$.

Consequently, by combining (\ref{equAS}), (\ref{eqwEq}), (\ref{eqmount}), (\ref{eqwGrad}), (\ref{eqz}), (\ref{eqwGrad1}), (\ref{eqomega}) and (\ref{eqwGrad2}), we infer that
\begin{equation}\label{equt}
u_t=o( e^{r\xi})\ \textrm{as} \ \xi \to +\infty,
\end{equation}
uniformly in $y\in \mathbb{R}^{N-1}$ and $t\in \mathbb{R}$.

\underline{\textbf{The sweeping argument.}} We are now in position to apply a sweeping argument in order to show the desired relation (\ref{eqdesire}).
Let \[\psi=u_t.\]
Differentiation of (\ref{equPDE}) with respect to $t$ yields
\begin{equation}\label{equPDEt}
  \psi_t=\psi_{\xi\xi}+c\psi_{\xi}+\Delta_y\psi+f'(u)\psi,\ \ \xi\in \mathbb{R}, \ y\in \mathbb{R}^{N-1},\ t\in \mathbb{R},
\end{equation}
Moreover, in terms of $\psi$ (\ref{equt}) becomes
\begin{equation}\label{equtpsi}
\psi=o( e^{r\xi})\ \textrm{as} \ \xi \to +\infty,\ \textrm{uniformly in}\ y\in \mathbb{R}^{N-1}\ \textrm{and} \ t\in \mathbb{R}.
\end{equation}

We observe that since $e^{r\xi}$ is a solution of (\ref{eqSupreme}) (keep in mind (\ref{eqr})), and $f$ satisfies (\ref{eqarxidia}), it is a supersolution of (\ref{equPDEt}) (the fact that it is a supersolution of (\ref{equPDE}) is well known even under the weaker condition (\ref{eqKPP})).
Indeed,
we have
\begin{equation}\label{eqsuperaras}
  (e^{r\xi})_t-(e^{r\xi})_{\xi\xi}-c(e^{r\xi})_{\xi}-\Delta_y e^{r\xi}-f'(u)e^{r\xi}=\left(f'(0)-f'(u)\right)e^{r\xi}\geq 0.
\end{equation}

Armed with the above information, we will show that $\psi\equiv 0$ by adapting Serrin's sweeping principle  (see \cite[Thm. 2.7.1]{sat} for the elliptic case).
Let us consider the set
\[\Lambda = \left\{
\lambda \geq 0 \  : \
\mu e^{r\xi} \geq \psi \ \textrm{in}\ \mathbb{R}^{N+1}
\ \textrm{for every}\ \mu \geq \lambda
\right\}.
\]	
Our goal is to show that $\Lambda = [0, \infty)$, which will yield $\psi \leq 0$. We can
also apply the same argument, with $\psi$ replaced by $-\psi$, to obtain $\psi \geq 0 $ and therefore
conclude.

We first prove that $\Lambda \neq \emptyset$, and thus by continuity
\begin{equation}\label{eqtildia}
  \Lambda = [\tilde{\lambda}, \infty)\ \textrm{for some}\ \tilde{\lambda}\geq 0.
\end{equation}
To this end, we note that since
$u$ is a bounded solution of (\ref{equPDE}) ($0<u<1$ as a matter of fact), standard interior estimates for linear
parabolic equations \cite{lady,lieber} and Sobolev embeddings imply that
\begin{equation}\label{eqOGholder}u\  \textrm{is bounded in}\
C^{2+\theta,1+\theta/2}
(\mathbb{R}^{N} \times \mathbb{R})\ \textrm{for any}\ \theta \in (0, 1).\end{equation}
So, this implies that
\begin{equation}\label{eqLinfty}
  \psi\in
L^\infty
(\mathbb{R}^{N} \times \mathbb{R}).
\end{equation}
The above relation and (\ref{equtpsi}) yield that there exists a  $\bar{\lambda}\gg 1$ such that
\[
\bar{\lambda} e^{r\xi} \geq \psi \ \textrm{in}\ \mathbb{R}^{N+1}.
\]
Hence, relation (\ref{eqtildia}) holds for some $\tilde{\lambda}\in [0,\bar{\lambda}]$.
For future reference, we note  that
\begin{equation}\label{eqfuturas}
  \psi \leq \tilde{\lambda} e^{r\xi}  \ \textrm{in}\ \mathbb{R}^{N+1}.
\end{equation}

In order to establish that $\tilde{\lambda}=0$, as desired, we will argue by contradiction. So,
let us suppose that $\tilde{\lambda}>0$. To show that this is absurd, by the definition of the
set $\Lambda$ and (\ref{eqtildia}) it suffices to prove that there exists a small $ \delta \in (0,\tilde{\lambda}/2)$ such that
\[
  (\tilde{\lambda} - \delta)e^{r\xi} > \psi(\xi,y, t),\ (\xi,y, t) \in \mathbb{R}^{N+1}.
\]
Suppose the above relation were false. Then, we could find $\lambda_n < \tilde{\lambda} $ with $\lambda_n\to \tilde{\lambda}$˜,
$\xi_n\in \mathbb{R}$, $y_n\in \mathbb{R}^{N-1}$ and  $t_n\in \mathbb{R}$  such that
\begin{equation}\label{eq35}\psi(\xi_n,y_n, t_n) \geq \lambda_ne^{r\xi_n},\ n \geq 1. \end{equation}
By virtue of (\ref{equtpsi}), (\ref{eqLinfty}), and our assumption that $\tilde{\lambda} >˜0$, we infer that the sequence
$\{\xi_n\}$ is bounded. Hence, passing to a subsequence if necessary, we may assume
that\begin{equation}\label{eqxn}
      \xi_n \to \xi_\infty \in \mathbb{R}.
    \end{equation}

Let us
now consider the  translated functions
\[U_n(\xi,y, t) = u(\xi,y+y_n, t + t_n)\ \textrm{and}\ \Psi_n(\xi,y, t)  = \psi(\xi,y+y_n, t + t_n),\ n \geq 1.\]
Clearly, $0< U_n< 1$   satisfies (\ref{eqOGholder}) uniformly with
respect to $n$; while  $\Psi_n$ solves
\begin{equation}\label{eqpsora}\Psi_t = \Psi_{\xi\xi}+c\Psi_{\xi}+\Delta_y\Psi+f'(U_n)\Psi,\ \xi\in \mathbb{R}, \ y\in \mathbb{R}^{N-1},\ t\in \mathbb{R},\end{equation}
and \begin{equation}\label{eqpsora2}
\Psi_n \     \textrm{is uniformly bounded with respect to}\ n\ (\textrm{recall (\ref{eqLinfty})}).
    \end{equation} We also note that (\ref{equtpsi}) becomes
\[
\Psi_n=o( e^{r\xi})\ \textrm{as} \ \xi \to +\infty,\ \textrm{uniformly in}\ y\in \mathbb{R}^{N-1}, \ t\in \mathbb{R}\ \textrm{and}\ n\geq 1.
\]
Moreover, from (\ref{eq35}) and (\ref{eqfuturas}) we obtain
\begin{equation}\label{eq38}
  \Psi_n(\xi_n,0, 0) \geq  \lambda_ne^{r\xi_n}\ \textrm{and}\ \Psi_n(\xi,y, t) \leq \tilde{\lambda}e^{r\xi},\ \forall\ (\xi,y, t) \in \mathbb{R}^{N+1},\ n\geq 1,
\end{equation}
respectively.

By the aforementioned uniform H\"{o}lder estimates for $U_n$ and a standard diagonal-compactness argument, passing to a further subsequence if necessary, we may
assume that
$U_n \to U_\infty$ in $C^{2,1}
_{loc}(\mathbb{R}^{N+1})$
for some $0\leq U_\infty\leq 1$  (actually, $U_\infty$ solves (\ref{equPDE}) and satisfies (\ref{equAS})-(\ref{eqw}) but we will not need this information).
In turn, by (\ref{eqpsora}), (\ref{eqpsora2}) and standard parabolic estimates, we deduce that $\Psi_n \to \Psi_\infty$ in $C^{2,1}
_{loc}(\mathbb{R}^{N+1})$
for some $  \Psi_\infty \in L^\infty(\mathbb{R}^{N+1})$ that solves
\begin{equation}\label{eqcontraduEq}
  \Psi_t = \Psi_{\xi\xi}+c\Psi_{\xi}+\Delta_y\Psi+f'(U_\infty)\Psi,\ \xi\in \mathbb{R}, \ y\in \mathbb{R}^{N-1},\ t\in \mathbb{R},
\end{equation}
and satisfies
\begin{equation}\label{eqcontradia}
\Psi_\infty=o( e^{r\xi})\ \textrm{as} \ \xi \to +\infty,\ \textrm{uniformly in}\ y\in \mathbb{R}^{N-1}, \ t\in \mathbb{R}.
\end{equation}
Moreover, recalling (\ref{eqxn}), by letting $n\to \infty$ in (\ref{eq38}) we get
\[\Psi_\infty(\xi_\infty,0, 0) \geq  \tilde{\lambda}e^{r\xi_\infty}\ \textrm{and}\ \Psi_\infty(\xi,y, t) \leq \tilde{\lambda}e^{r\xi},\ (\xi,y, t) \in \mathbb{R}^{N+1}.\]
On the other hand, since $\tilde{\lambda}e^{r\xi}$ is a supersolution of (\ref{eqcontraduEq}) (by the same calculation as in (\ref{eqsuperaras}) but with $U_\infty$ in place of $u$), we infer by the strong maximum principle (see for instance \cite{lady,lieber}) that $\Psi_\infty \equiv \tilde{\lambda}e^{r\xi}$. However, the last relation  contradicts (\ref{eqcontradia}), and thus the proof of the theorem is complete.
\end{proof}

\textbf{Acknowledgments.} The author would like to thank IACM of FORTH, where
this paper was written, for the hospitality. This work has received funding from
the Hellenic Foundation for Research and Innovation (HFRI) and the General
Secretariat for Research and Technology (GSRT), under grant agreement No
1889.

\end{document}